\documentclass[]{iac}

\usepackage{times}  
\usepackage{helvet}  
\usepackage{courier}  
\usepackage[hyphens]{url}  
\usepackage{graphicx} 
\usepackage{cite}
\usepackage{caption} 
\DeclareCaptionStyle{ruled}{labelfont=normalfont,labelsep=colon,strut=off} 
\usepackage{algorithm}
\usepackage{algorithmic}

\usepackage{times}
\usepackage{epsfig}
\usepackage{graphicx}
\usepackage{amsmath}
\usepackage{amssymb}
\usepackage{array}
\usepackage[hidelinks]{hyperref}

\usepackage{newfloat}
\usepackage{listings}
\newcolumntype{L}[1]{>{\raggedright\arraybackslash}m{#1}}
\newcolumntype{C}[1]{>{\centering\arraybackslash}m{#1}}
\newcolumntype{R}[1]{>{\raggedleft\arraybackslash}m{#1}}

\begin{document}

\IACpaperyear{22}
\IACpapernumber{C1.LBA}  
\IACconference{73}
\IAClocation{Paris, France}
\IACdate{18-22 September 2022}
\IACcopyrightA{2022}{the authors}

\title{A Lambert's Problem Solution via the Koopman Operator with Orthogonal Polynomials}

\IACauthor{Julia Pasiecznik}{1}{}
\IACauthor{Simone Servadio}{2}{}
\IACauthor{Richard Linares}{3}{}

\IACaffiliation{SM Candidate, Department of Aeronautics and Astronautics, Massachusetts Institute of Technology, MA 02139, USA.}{jpasiecz@mit.edu}{1}
\IACaffiliation{Postdoctoral Associate, Department of Aeronautics and Astronautics, Massachusetts Institute of Technology, MA 02139, USA.}{simoserv@mit.edu}{2}
\IACaffiliation{Rockwell International Career Development Professor, Department of Aeronautics and Astronautics, Massachusetts Institute of Technology, MA 02139, USA.}{linaresr@mit.edu}{3}

\abstract{Lambert’s problem has been long studied in the context of space operations; its solution enables accurate orbit determination and spacecraft guidance. This work offers an analytical solution to Lambert’s problem using the Koopman Operator (KO). In contrast to previous methods in the literature, the KO provides the analysis of a nonlinear system by seeking a transformation that embeds the nonlinear dynamics into a global linear representation. Our new methodology to solve for Lambert solutions considers the position of the system's eigenvalues on the phase plane, evaluating accurate state transition polynomial maps for a computationally efficient propagation of the dynamics. 
The methodology used and multiple-revolution solutions found are
compared in accuracy and performance with other techniques found in
the literature, highlighting the benefits of the newly developed analytical approach over classical numerical methodologies.}{Lambert's Problem, Koopman Operator}

\maketitle



\section{Introduction}
Solving Lambert's problem concerns finding an optimal transfer orbit given an initial and final position vector and a time of flight \cite{lancaster}. It is an orbital boundary value problem that has been studied by mathematicians for years and is at the core of many astrodynamical problems. Finding solutions to Lambert's problem is key to modern spacecraft operations as these solutions can be used for orbit determination, multi-planet trajectories, targeting and rendez-vous. 

Various approaches have been designed to solve Lambert's problem as reviewed by de la Torre Sangr\`a and Fantino \cite{sangra}. The fundamental work was developed by Euler, Lambert, Lagrange, and Gauss \cite{bopp} \cite{lagrange} \cite{gauss}. Lancaster and Blanchard expanded on this work and solved Lambert's problem using Halley's cubic iteration process \cite{lancaster}. Gooding built upon Lancaster and Blanchard's method and developed a method that solved Lambert's problem with high precision within three iterations \cite{gooding}. Izzo also built upon Lancaster and Blanchard's method and developed a solver with a reduction in the overall solver complexity \cite{izzo}.
Izzo's Householder algorithm has been found to have the best ratio between speed, robustness and accuracy \cite{sangra}. 
The multiple-revolutions Lambert's problem includes optimal transfer orbits that have a number of complete revolutions before arriving at the final position vector. Furthermore, solutions to solving the multiple-revolutions Lambert's problem generally fall into two categories, 1) methods that use direct geometry and 2) methods that use universal variables \cite{arora}. Many of the approaches mentioned above are capable of solving multiple-revolutions Lambert's problem in addition to the single revolution case.

In this paper we develop a new analytical method to solve both single and multiple-revolutions Lambert's problem using the Koopman Operator. The KO provides the analysis of a nonlinear system by seeking a transformation that embeds the nonlinear dynamics into a global linear representation. The resulting solution is represented as a linear combination of the system’s eigenfunctions, since the KO obtains the spectral behavior of the system, in terms of its modes, frequencies, and eigenvalues. Thanks to this approach, this work provides a new methodology to solve for multiple-revolution Lambert solutions, by considering the position of the eigenvalues on the phase plane. 
The KO has been applied to celestial mechanics applications such as the zonal harmonics problem around a celestial body and the circular restricted three body problem \cite{arnaslinares}. It has been shown to find analytical solutions that are stable and accurate \cite{threebody}. 

The KO approach evaluates accurate state transition matrices (STM) for a computationally efficient propagation of the dynamics. Considering the selection of specific energy as the cost function, Lambert’s problem is first solved in this work without perturbation effects. A solution that considers the oblateness of the Earth ($J_2$ effects) is then evaluated by merging perturbation theory into the KO framework. The perturbed frequencies and perturbed eigenfunctions are found using higher order corrections of the solution to the unperturbed problem. As such, the perturbed solution is found by merging the unperturbed solution (solvable part) with highly nonlinear terms thanks to the use of perturbation theory. Using the unperturbed solution reduces the complexity involved in solving the perturbed problem by eliminating the need to perform the eigendecomposition of the perturbing matrix. The KO perturbation method has been shown to provide solutions comparable in accuracy to other perturbation methods such as the Poincar\'e Lindstedt method when applied to the zonal harmonics problem \cite{arnaslinares}.
This work utilizes a new set of orbital elements to obtain the KO eigendecomposition of the Lambert’s problem, where the dynamics are represented in a polynomial manner and projected onto a well-defined set of orthogonal basis functions. Thus, different transfer orbits are obtained as a linear combination of the KO eigenfunctions. The solutions found using this new method to solve Lambert's problem are compared in accuracy with other techniques found in the literature, highlighting the benefits of the newly developed analytical approach over classical numerical methodologies.

\section{Koopman Operator Theory}
The KO is an infinite dimensional linear operator that evolves observables of functions of states. We present a summary of KO theory and the methodology used to apply it to Lambert's problem accounting for $J_2$ effects.
We consider a classical nonlinear dynamical autonomous system given by the set of ordinary differential equations, 
\begin{equation}
\left\{\begin{array}{l}
\frac{d}{d t} \mathbf{x}(t)=\mathbf{f}(\mathbf{x}) \\
\mathbf{x}\left(t_{0}\right)=\mathbf{x}_{\mathbf{0}}
\end{array}\right.
\end{equation}
where $\mathbf{x} \in \mathbb{R}^{d}$ is the set of variables that describe the system evolving in time, $d$ is the number of dimensions of the space where the problem is defined, and $\mathbf{x}_{\mathbf{0}}$ represents the initial conditions at the initial time, $t_0$.
Let us consider an observable function $\boldsymbol{g}(\mathbf{x})$ of this system. The evolution of $g(\mathbf{x})$ is given by applying the KO, $\mathcal{K}$, onto the observable \cite{arnaslinares},
\begin{equation}
\begin{array}{c c}
\frac{d}{d t} \boldsymbol{g}(\mathbf{x})=\mathcal{K}(\boldsymbol{g}(\mathbf{x}))=\left(\nabla_{\mathbf{x}} \boldsymbol{g}(\mathbf{x})\right) \mathbf{f}(\mathbf{x})
\end{array}.
\end{equation}
This relation can be generalized to any observable of the dynamical system, 
\begin{equation}
\mathcal{K}(\cdot)=\left(\nabla_{\mathbf{x}} \cdot\right) \mathbf{f}(\mathbf{x}).
\end{equation}
The KO is an infinite dimensional linear operator and hence the vector space of the observable functions given by $\mathcal{F}, \text { where } g_{i}(\mathbf{x}) \in \mathcal{F}$ and $i \in\{1, \ldots, d\}$, is an infinite dimensional Hilbert space. 
Since we can not work with an infinite Hilbert space, we limit the number of dimensions and evolve the system using a finite set of basis functions spanning a finite subspace \cite{brunton2016}. This extended space is still able to capture the non-linearities of the original system however as an approximation. 
Hence, the KO provides an approximate linear description of nonlinear dynamical systems, allowing for efficient computation of solutions to problems in nonlinear astrodynamical systems.

Any set of eigenfunctions of the KO, $\boldsymbol{\Phi}(\mathbf{x})  \in \mathcal{F}$, provides a Koopman invariant finite subspace, $\mathcal{F}_D$, with dimension $m$, in which to evolve the system \cite{linareskoopman},
\begin{equation}
\label{eigenfunctionevolution}
\mathcal{K}(\boldsymbol{\Phi}(\mathbf{x}))=\frac{\mathrm{d}}{\mathrm{d} t} \boldsymbol{\Phi}(\mathbf{x})=\boldsymbol{\Lambda} \boldsymbol{\Phi}(\mathbf{x})
\end{equation}
where $\boldsymbol{\Lambda}$ is a diagonal matrix containing the eigenvalues $\boldsymbol{\Lambda}=\operatorname{diag}\left(\left[\lambda_1, \ldots, \lambda_m\right]\right)$. 
The evolution of the eigenfunctions in time is given by the solution to equation \eqref{eigenfunctionevolution},
\begin{equation}
\boldsymbol{\Phi}(t)=\exp (\boldsymbol{\Lambda} t) \boldsymbol{\Phi}\left(t_{0}\right).
\label{evolutioneigenfunction}
\end{equation}
Once the eigenfunctions of the system are found, equation \eqref{evolutioneigenfunction} can be used to solve the complete system. This approximate linear description provided by the KO allows us to take advantage of the many mathematical tools provided by linear algebra, such as computing the state transition matrices. However, before computing the full solution by making use of such tools, we first compute the set of eigenfunctions of the KO using the Galerkin method. 
\subsection{Computation of the Koopman Eigenfunctions}
\label{sec:galerkinmethod}
The Galerkin method provides the approximated linearized system of equations needed to perform the linearization in the KO theory. It projects the system onto a set of orthogonal functions that span the finite subspace of the KO theory. The method has been previously introduced by Arnas and Linares \cite{arnaslinares}, we include a summary of key parts of the formulation here for its application to Lambert's problem.
Let the dynamical system be represented by $f(\mathbf{x})$ and let $h(\mathbf{x})$ represent the basis functions defined in a domain $\Omega$. The inner product between these two functions provides the projection,
\begin{equation}
\langle f, h\rangle=\int_{\Omega} f(\mathbf{x}) h(\mathbf{x}) w(\mathbf{x}) d \mathbf{x}
\end{equation}
where $w(\mathbf{x})$ is a positive weighting function associated with the basis functions.
We use orthonormal Legendre polynomials as the set of orthonormal basis functions due to the advantages they give in simplifying the inner product computation. The Legendre polynomials are defined in a bounded domain $([-1,1])$ and the weighting function associated with the Legendre polynomials is a constant.
Let $\mathbf{L}(\mathbf{x})$ denote the set of all basis functions. Written in vector form, the set is given by $L(\mathbf{x})=\left[L_{1}(\mathbf{x}), \ldots, L_{m}(\mathbf{x})\right]^{T}$ where $L_i(\mathbf{x})$ is a multidimensional basis function of the set. The orthonormal basis functions can be given by the Legendre polynomials, by expressing $L_i(\mathbf{x})$ as,
\begin{equation}
L_{i}(\mathbf{x})=\prod_{j=1}^{m} P_{n_{j}}\left(\mathrm{x}_{j}\right)
\end{equation}
where $P_{n_{j}}\left(\mathrm{x}_{j}\right)$ denotes the Legendre polynomial of order $n_{j}$ acting on variable $\mathrm{x}_{j}$ \cite{arnaslinares}. 
The total order of the polynomial, c, is given by adding all of the individual orders of the Legendre polynomials unique to each basis function,
\begin{equation}
c = \sum_{j=1}^{m} n_{j}
\end{equation}
By applying the KO onto $L_{i}(\mathbf{x})$ we can obtain its total derivative,
\begin{equation}
\mathcal{K}\left(L_{i}(\mathbf{x})\right)=\frac{d}{d t} L_{i}(\mathbf{x})=\left(\nabla_{\mathbf{x}} L_{i}(\mathbf{x})\right) \mathbf{f}(\mathbf{x}).
\end{equation}
Taking advantage of the inner product operation, we can approximate the derivative of the basis functions using the Koopman matrix, $K$. The components of the matrix, $K_{ij}$ give the projection of the derivative of the $i$th basis function onto the $j$th basis function,
\begin{equation}
\begin{aligned}
K_{i j}&=\left\langle\left(\nabla_{\mathbf{x}} L_{i}(\mathbf{x})\right) \mathbf{f}(\mathbf{x}), L_{j}(\mathbf{x})\right\rangle\\
&=\int_{\Omega}\left(\nabla_{\mathbf{x}} L_{i}(\mathbf{x})\right) \mathbf{f}(\mathbf{x}) L_{j}(\mathbf{x}) w(\mathbf{x}) d \mathbf{x}.
\end{aligned}
\end{equation}
Using this formulation for all combinations of basis functions, we arrive at the best approximated linear system of the dynamics for the selected basis functions given by\cite{arnaslinares},
\begin{equation}
\frac{d}{d t} \mathbf{L}(\mathbf{x})=K \mathbf{L}(\mathbf{x}).
\end{equation}
Increasing the order of the basis functions allows the KO matrix to better encapsulate the non-linearities of the dynamical system and hence improve the accuracy of the solution at the cost of increasing the size of the KO matrix. However, the KO matrix needs only to be computed once to be applied to any initial condition.


In order to represent the observables in terms of the system's set of eigenfunctions, 
$\boldsymbol{\Phi}(\mathbf{x})$, we must use a transformation matrix, $B$, to make a projection onto the basis functions. In component form this transformation matrix is \cite{threebody},
\begin{equation}
b_{i \ell}=\left\langle g_i, \phi_{\ell}\right\rangle=\int_{\Omega} g_{i}(\mathbf{x}) \phi_{\ell}(\mathbf{x}) w(\mathbf{x}) d \mathbf{x}
\end{equation}
where $b_{i\ell}$ are the coefficients of the matrix $B$, which has size $q \times m$. Here $q$ represents the number of observables of the system. 
$B$ is known as the Koopman modes matrix and it relates the observables with the eigendecomposition of the dynamics.
The observables can now be written in terms of the Koopman modes matrix,
\begin{equation}
\boldsymbol{g}(\mathbf{x})=B \mathbf{\Phi(x)}.
\end{equation}
The evolution of the observables in time can be found using the evolution of the eigenfunctions of the system, 
\begin{equation}
\boldsymbol{g}(\mathbf{x}(t))=B \boldsymbol{\Phi}(\mathbf{x}(t))=B \exp (\boldsymbol{\Lambda} t) \boldsymbol{\Phi}\left(\mathbf{x}\left(t_0\right)\right).
\label{eqn:evolutioneigen}
\end{equation}
Let $P$ be the matrix containing all the eigenvectors of the Koopman matrix, $K$. The coefficient vectors of the KO eigenfunctions are the eigenvectors of the Koopman matrix with the same eigenvalues, 
\begin{equation}
P K=\Lambda P.
\end{equation}
We can express the result in equation \eqref{eqn:evolutioneigen} in terms of $P$ and the basis functions $\boldsymbol{L}$ \cite{threebody},
\begin{equation}
\boldsymbol{g}(\mathbf{x}(t))=B \boldsymbol{\Phi}(\mathbf{x}(t))=B \exp (\boldsymbol{\Lambda} t) P \boldsymbol{L}\left(\mathbf{x}\left(t_0\right)\right).
\label{eqn:observableevolution}
\end{equation}

Let $a_{i \ell}$ represent the projection of the observables onto the basis functions, $a_{i \ell}=\left\langle g_i, L_{\ell}\right\rangle$. We can use this projection to represent the observables in terms of the basis functions,
\begin{equation}
\boldsymbol{g}(\mathbf{x})=A \boldsymbol{L}(\mathbf{x}) \rightarrow g_i(\mathbf{x}) \approx \sum_{\ell=1}^m a_{i \ell} L_{\ell}(\mathbf{x})
\label{eqn:evolutionobs}
\end{equation}
where $A_{i \ell} = a_{i \ell}$. 
The solution to the set of differential equations given in equation \eqref{eqn:evolutionobs} is given by,
\begin{equation}
\boldsymbol{g}(\mathbf{x}(t))=A \boldsymbol{L}(\mathbf{x}(t))=A P^{-1} \exp (\boldsymbol{\Lambda} t) P \boldsymbol{L}\left(\mathbf{x}\left(t_0\right)\right).
\label{eqn:evolutionsolution}
\end{equation}
This gives the evolution of the observables over time and the representation of any polynomial function governed by the dynamics of the system. 
We can now use equation \eqref{eqn:observableevolution} to obtain the state transition matrix (STM) of the system.
Given the state vector at an initial time, i.e. the initial position vector, the STM allows us to calculate the general solution of the linearized dynamical system that the KO formulation provides.
In the case where the function governing the observable is given by the identity, $\boldsymbol{g}(\mathbf{x})=\mathbf{x}$, the STM is given by,
\begin{equation}
\frac{\partial \mathbf{x}\left(t_f\right)}{\partial \mathbf{x}\left(t_0\right)}=B \exp \left(\boldsymbol{\Lambda} t_f\right) P \frac{\partial \boldsymbol{L}\left(\mathbf{x}\left(t_0\right)\right)}{\partial \mathbf{x}\left(t_0\right)}
\label{eqn:STM}
\end{equation}
where $t_f$ is the final time of the dynamics. Now we will use the KO formulation to find solutions to Lambert's problem. Once again, the KO matrix only needs to be computed once for the dynamics of the system and then can be used to propagate any initial condition, hence making it easy to find the minimum energy orbit transfer of Lambert's problem for a given time range.


\section{Defining the Set of Variables used to Solve Lambert's Problem}
\label{Sec:Formulation to Solve Lambert's Problem}
In order for the KO theory to be well applied to a nonlinear system, the system must be close to a linear system with a small perturbation. This enables the limited number of basis functions to capture the non-linearities and for the KO theory to efficiently find the solution to the astrodynamical problem. A set of new orbital elements have been introduced in the zonal formulation work of Arnas and Linares \cite{arnaslinares} that is close to linear and can be used to solve Lambert's problem with and without $J_2$ perturbations. We consider a satellite orbiting Earth as the system in which we want to solve Lambert's problem.
The Hamiltonian of the system subject to gravitational perturbations given by $J_2$ effects is given in spherical coordinates as
\begin{equation}
\begin{aligned}
\mathcal{H}=&\frac{1}{2}\left(p_r^2+\frac{p_{\varphi}^2}{r^2}+\frac{p_\lambda^2}{r^2 \cos ^2(\varphi)}\right)\\
&-\frac{\mu}{r}+\frac{1}{2} \mu R_{\oplus}^2 J_2 \frac{1}{r^3}\left(3 \sin ^2(\varphi)-1\right)
\end{aligned}
\end{equation}
where $r$ is the radial distance from the orbiting satellite to the center of the Earth, $\mu$ represents the gravitational constant of the Earth, $R_{\oplus}$ represents the Earth's equatorial radius, $\varphi$ is the latitude of the satellite, $\lambda$ is the longitude of the satellite, and $J_2$ is the zonal term of the Earth's gravitational potential.
Furthermore, $p_r$, $p_{\varphi}$, and $p_\lambda$ represent the conjugate momenta of the system and are given by
\begin{equation}
\begin{aligned}
p_r &=\dot{r} \\
p_{\varphi} &=r^2 \dot{\varphi} \\
p_\lambda &=r^2 \cos ^2(\varphi) \dot{\lambda}
\end{aligned}
\end{equation}
where $\dot{x}$ is the derivative of the variable $x$ with respect to time. 

The set of new orbital elements defined by Arnas and Linares \cite{arnaslinares} are defined by the following transformations,
\begin{equation}
\begin{aligned}
\Lambda &=\sqrt{\frac{R_{\oplus}}{\mu}}\left(\frac{p_\theta}{r}-\frac{\mu}{p_\theta}\right)\\
\eta &=p_r \sqrt{\frac{R_{\oplus}}{\mu}} \\
s &=\sin (\varphi) \\
\gamma &=\frac{p_{\varphi}}{p_\theta} \cos (\varphi)  \\
\kappa &=\sqrt{\frac{\mu R_{\oplus}}{p_{\varphi}^2+\frac{p_\lambda^2}{\cos ^2(\varphi)}}} \\
\beta &=\lambda-\arcsin \left(\tan (\varphi) \sqrt{\frac{p_\lambda^2}{p_\theta^2-p_\lambda^2}}\right)  \\
\chi &=\frac{p_\lambda}{p_\theta^4} \frac{\left(\mu R_{\oplus}\right)^{3 / 2} }{\sin ^2(\varphi)+\frac{p_\varphi^2}{p_\theta^2} \cos ^2(\varphi)} \\
\rho &=\frac{p_\lambda}{p_\theta}
\end{aligned}
\label{eq:orbitalelements}
\end{equation}
where 
\begin{equation}
p_\theta=\sqrt{p_{\varphi}^2+\frac{p_\lambda^2}{\cos ^2(\varphi)}}
\end{equation}
is the angular momentum of the orbit.

The corresponding Hamilton equations given in the set of new orbital elements and taking into account $J_2$ effects are 
\begin{equation}
\begin{aligned}
\frac{d \Lambda}{d \theta} &=-\eta-3 J_2 s \gamma \kappa^3(\Lambda+\kappa)(\Lambda+2 \kappa) \\
\frac{d \eta}{d \theta} &=\Lambda+\frac{3}{2} J_2 \kappa^3(\Lambda+\kappa)^2\left(3 s^2-1\right) \\
\frac{d s}{d \theta} &=\gamma \\
\frac{d \gamma}{d \theta} &=-s-3 J_2 s \rho^2 \kappa^3(\Lambda+\kappa) \\
\frac{d \kappa}{d \theta} &=3 J_2 s \gamma \kappa^4(\Lambda+\kappa) \\
\frac{d \beta}{d \theta} &=-3 J_2 s^2 \chi(\Lambda+\kappa) \\
\frac{d \chi}{d \theta} &=12 J_2 s \gamma \chi \kappa^3(\Lambda+\kappa)+6 J_2 s \rho \chi^2(\Lambda+\kappa) \\
\frac{d \rho}{d \theta} &=3 J_2 s \gamma \rho \kappa^3(\Lambda+\kappa)
\end{aligned}
\label{eqn:elemwithj2}
\end{equation}
The linear contributions to these variables that do not have any terms proportional to $J_2$ give the system of equations for the classical two body problem without $J_2$ effects.
The reader is referred to the zonal formulation made by Arnas and Linares in \cite{arnaslinares} for a detailed derivation of these variables. 



To make use of the orbital elements presented in equation \eqref{eq:orbitalelements}, we make a transformation of variables from Cartesian coordinates, to spherical coordinates, to the set of new orbital elements. Once the solution is found using the KO in the new set of orbital elements, we use the inverse of the transformations to get back to cartesian coordinates. These transformations allow us to formulate the problem statements in familiar cartesian or spherical coordinates.
Now let us solve Lambert's problem in the single revolution case.


\section{Single Revolution Lambert's Problem}
\subsection{Classic Two Body Dynamics}
\label{sec:withoutJ2}
The solution to the single revolution Lambert's problem concerns finding a direct optimal orbit transfer from two position vectors and a given time of flight. We first solve the problem without taking into account the effects of $J_2$ perturbations from Earth. We present the problem and its solution in Cartesian coordinates.

We consider a satellite orbiting Earth with initial conditions taken from Example 5.2 of the Orbital Mechanics for Engineering Students textbook by H. Curtis \cite{curtis}.
The initial and desired final position vectors are $\mathbf{r_0} = [5000, \hspace{0.2cm} 10,000, \hspace{0.2cm}2100]$, and $\mathbf{r_f} = [-14,600,\hspace{0.2cm} 2500,\hspace{0.2cm} 7000]$ respectively. 
The time of flight is $3600$ seconds.
The angle between the initial and final position vectors, $\Delta{\theta}$, is given by:
\begin{equation}
    \Delta{\theta} = \arccos{\frac{\mathbf{r_0}*\mathbf{r_f}}{|{\mathbf{r_0}}||{\mathbf{r_f}}|}},
    \label{eq:deltatheta}
\end{equation}
$\Delta{\theta}$ is also known as the true anomaly deviation. 

After transforming the given initial state into the new set of orbital elements described in Section \ref{Sec:Formulation to Solve Lambert's Problem}, the propagation is performed in the Koopman framework according to equation \eqref{eqn:evolutionsolution}, where the time of flight has been substituted.
Using a function that makes use of the Levenberg-Marquardt algorithm (LM) and an initial circular velocity guess, the initial velocity, $\mathbf{v_0}$, that yields the desired final position within the time constraint specified is solved for. 
Passing different time constraints through the LM function yields different initial velocities. However, passing in different time constraints does not require recalculating the KO matrix since it is calculated only once by projecting the dynamics of the system onto the basis functions.

\begin{figure}
\begin{center}
\includegraphics[width=0.9\linewidth]{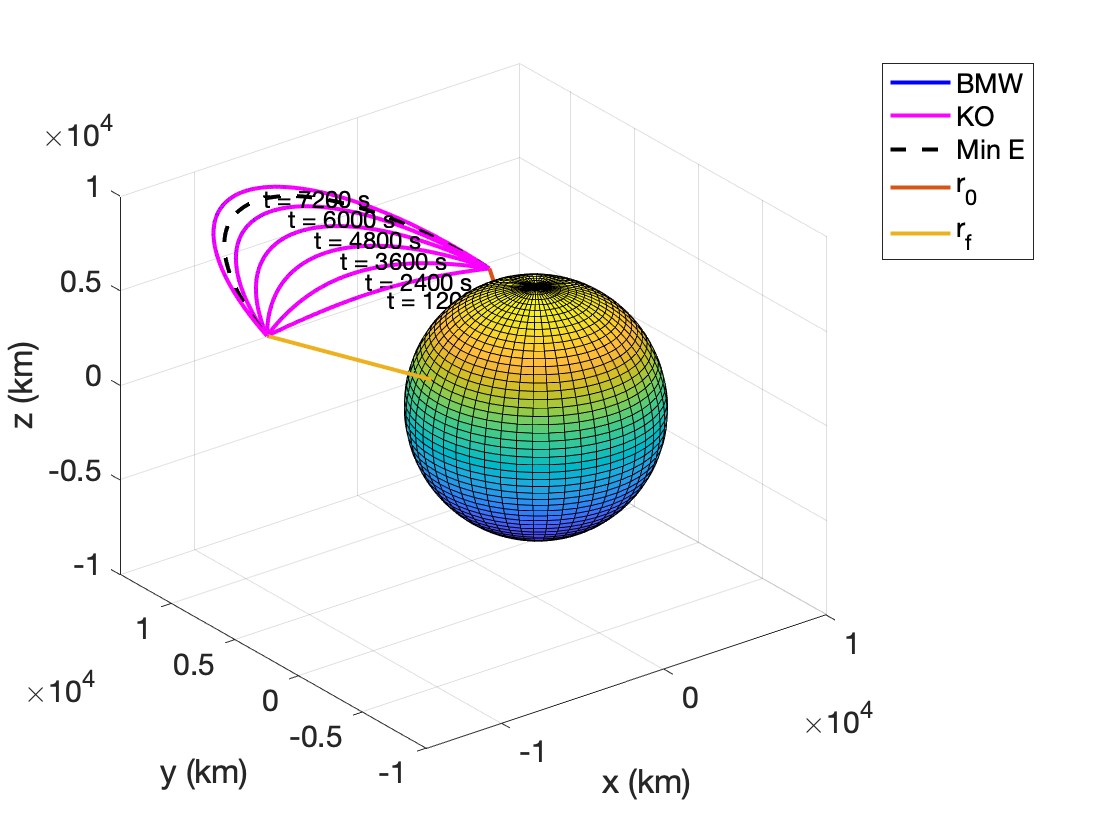}
\caption{Single Revolution Optimal Trajectories without $J_2$ Effects}
\label{fig:SingleRevolution}
\end{center}
\end{figure}

For our example, we plot the resulting orbital trajectories found for six different times, equally spaced between $1200$s and $7200$s in Figure \ref{fig:SingleRevolution}.
The two body problem in Cartesian coordinates, solved using the algorithm by Bate, Mueller and White \cite{bate} (reported as BMW in the Figures) and Bond and Allman \cite{bond} is plotted for comparison. 

\begin{figure}
\begin{center}
\includegraphics[width=0.9\linewidth]{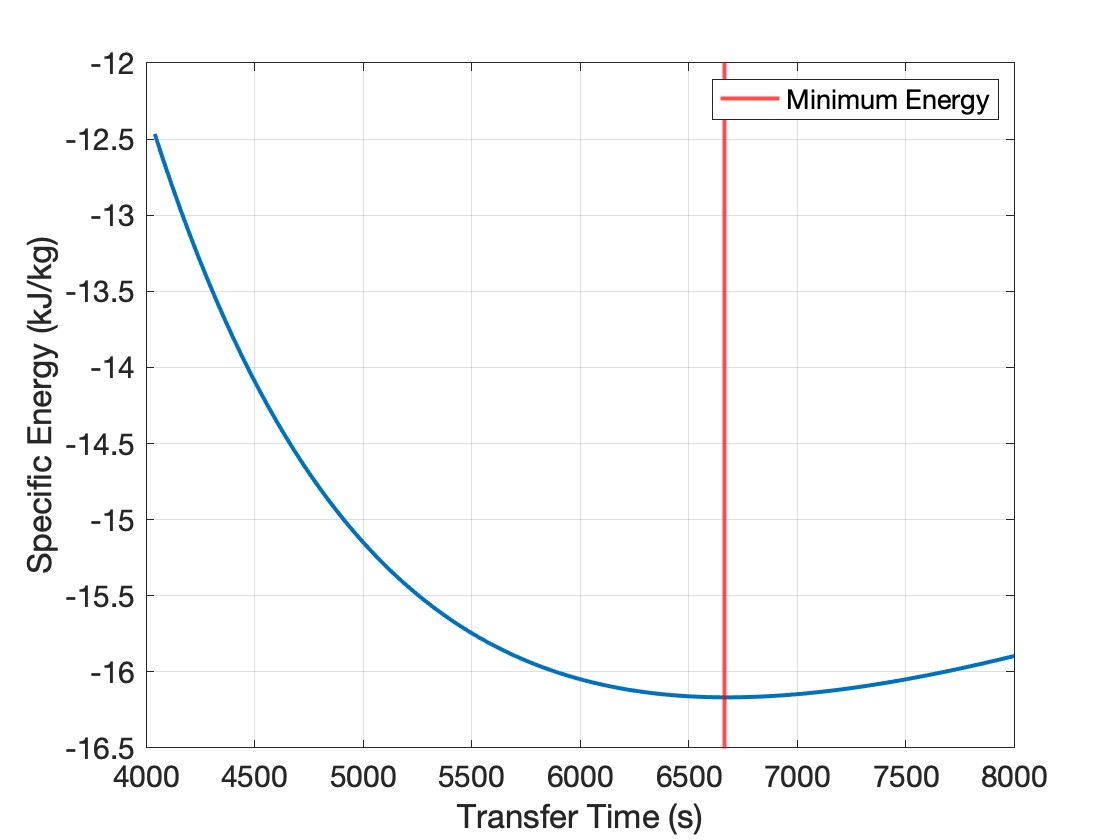}
\caption{Specific Energy of Each Transfer Orbit}
\label{fig:minenergy}
\end{center}
\end{figure}

We select specific energy as the cost function to find the transfer orbit with the minimum specific energy and hence solve Lambert's problem. 
The specific energy is given by
\begin{equation}
E=-\frac{\mu}{2 a}
\end{equation}
where $a$ is the semi major axis and is calculated for each transfer orbit from the initial cartesian coordinates
using 
\begin{equation}
a = \frac{1}{\frac{2}{|{\mathbf{r_0}}|} - \frac{|{\mathbf{v_0}}|^2}{\mu}}.
\end{equation}
Negative specific energy solutions indicate elliptical transfer orbits. Finding the specific energy of each possible transfer orbit allows us to find the minimal specific energy transfer orbit. We find the minimum specific energy transfer orbit corresponds to a time of flight of $6630$s which corresponds to a semi-major axis of $12,328$ km and a minimum specific energy of $-16.2$    kJ/kg. We plot the specific energies for transfer times near the time corresponding to the minimum specific energy solution in Figure \ref{fig:minenergy}.
The minimum specific energy transfer orbit is plotted in Figure \ref{fig:SingleRevolution}.

\subsection{With $J_2$ Perturbations}
\begin{figure}
\begin{center}
\includegraphics[width=0.9\linewidth]{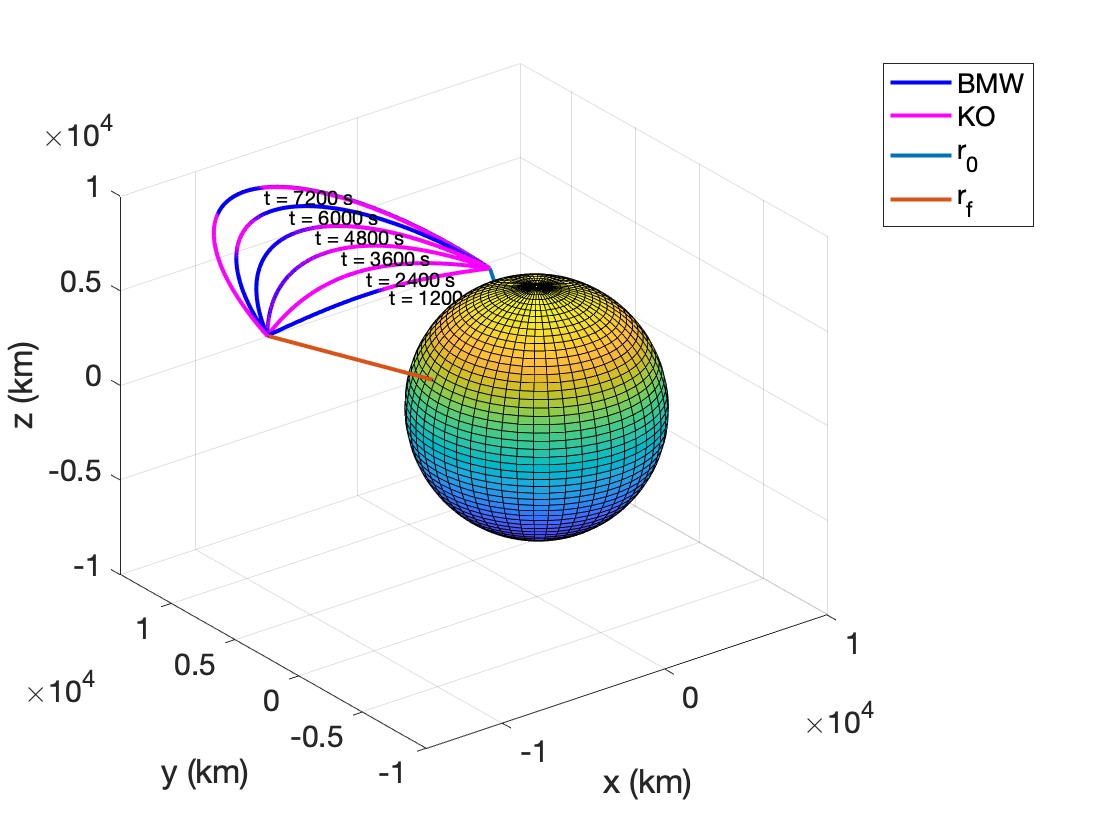}
\caption{Single Revolution Optimal Trajectories with $J_2$ Effects}
\label{fig:SingleRevolutionJ2}
\end{center}
\end{figure}

Perturbations due to the oblateness of the Earth are taken into account by including the terms proportional to $J_2$ in the equations that represent the change in the orbital elements given by \eqref{eqn:elemwithj2}. In order to calculate the optimal trajectories, we increase the order of the basis functions of the Koopman operator up to order 7. As the maximum order of the monomials in the $J_2$ dynamics is 7, having order 7 of the Legendre polynomials guarantees to cover any nonlinear contributions. Increasing the order improves the accuracy of the solution by being able to better take into account the non-linearities of the system. We plot the resulting orbital trajectories found for the same six different times, between $1200$s and $7200$s in Figure \ref{fig:SingleRevolutionJ2}.

\begin{figure}
\begin{center}
\includegraphics[width=0.9\linewidth]{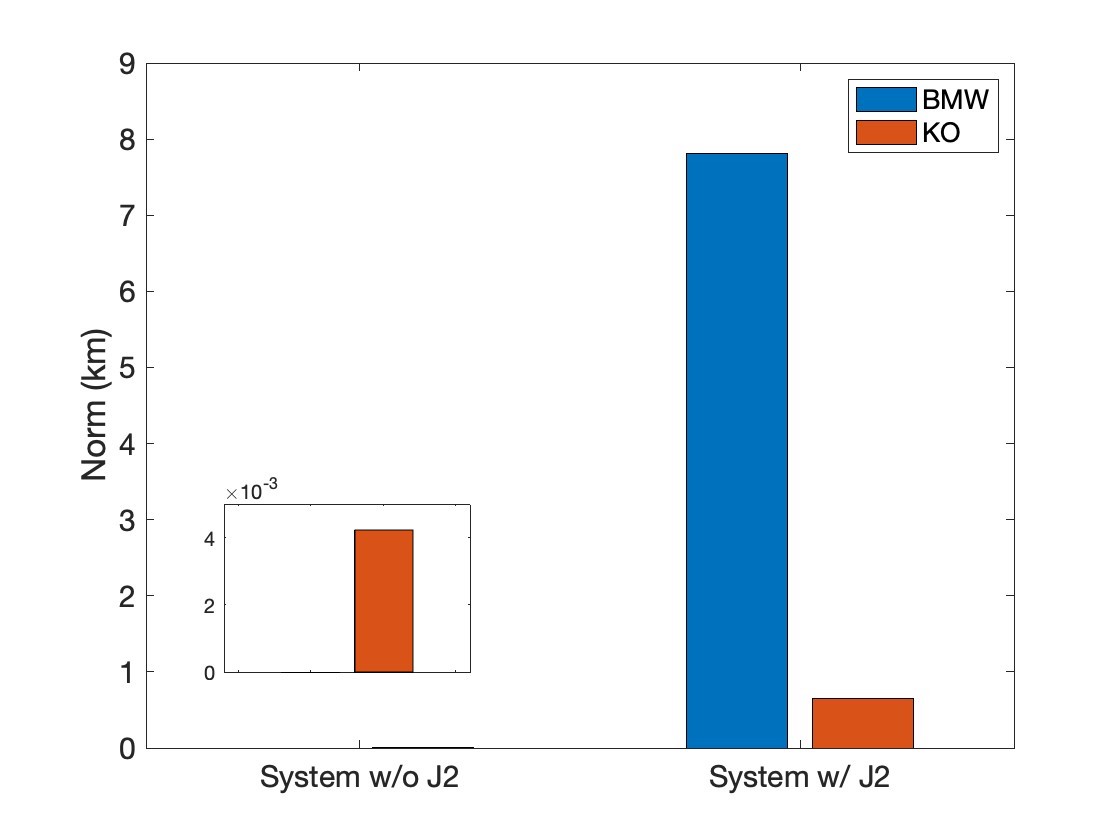}
\caption{Norm of difference Between Attained Final Position and Desired Final Position for the Propagation of Velocities in the Dynamical System with and without $J_2$  effects.}
\label{fig:comparisonplotv2}
\end{center}
\end{figure}

We propagate the initial velocities found using both BMW and KO that give the transfer orbits for the single revolution case without taking into account $J_2$ effects for the time constraint of $t = 3600$s. 
The propagator takes the initial state of the system and propagates it in time to find the final state after the given time of flight to validate results. 
The propagation is done in a two body system without $J_2$ effects using a variable order Adams-Bashforth-Moulton PECE solver (ABM).
We find the norm of the difference between the final position that the satellite reaches and the desired final position, $\mathbf{r_f}$, for the velocity found using the methodology by BMW to be $1.1997 \times 10^{-11}$ km and for the velocity found using KO to be $0.00421$ km.
Subsequently, the initial velocity found using BMW for the single revolution case is propagated in a two body system with $J_2$ effects using ABM. 
The norm of the difference between the final position that the satellite reaches and the desired final position is $7.81$ km.
We also propagate the initial velocity found using the KO for the single revolution case where we have taken into account $J_2$ effects in a two body system with $J_2$ effects using ABM. In this case, the norm of the difference between the final position that the satellite reaches and the desired final position is $0.655$ km.
Hence, the velocity found using the KO for the single revolution case propagated in a two body system with $J_2$ effects gets the satellite closer to the desired final position than using the velocity found using the BMW algorithm which does not take into account $J_2$ effects.
We plot the results comparing the norms in Figure \ref{fig:comparisonplotv2}.


\section{Multiple Revolutions Lambert's Problem}
\subsection{Without $J_2$ Perturbations}
\begin{figure}
\begin{center}
\includegraphics[width=0.9\linewidth]{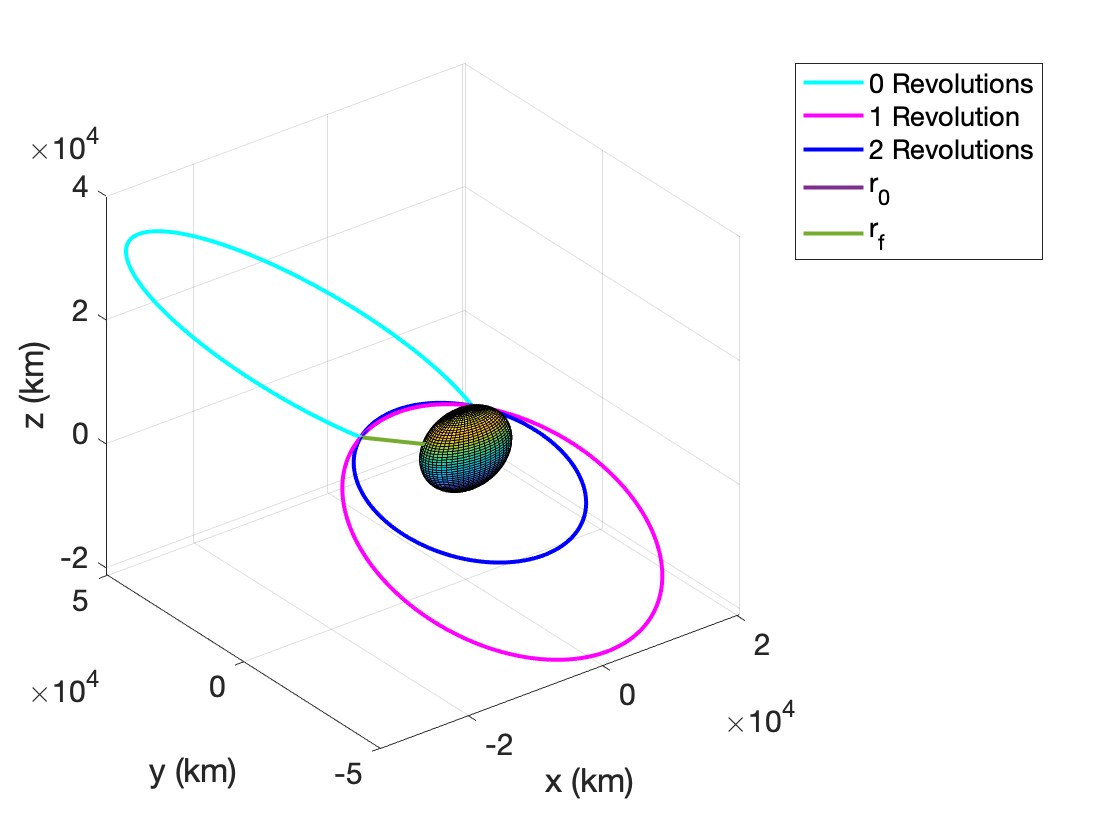}
\caption{Multiple Revolutions Orbit Trajectories without $J_2$ Effects for $t = 59,952$s}
\label{fig:MultipleRevolution}
\end{center}
\end{figure}

Solving Lambert's Problem using multiple revolutions involves finding optimal orbit transfers that complete more than one revolution around the gravitational body when travelling from the initial position vector to the final position vector within a given amount of time.
We consider the case from Section \ref{sec:withoutJ2} without taking into account effects from $J_2$ perturbations. 
The period of such an orbit travelling from $\mathbf{r_0}$ to $\mathbf{r_f}$ is calculated to be $T = 28,175$s. We choose a large transfer time corresponding to two orbital periods $59,952$s, in order to solve for multiple revolutions.
These multiple revolutions solutions are found by adding $2 \pi N$ to $\Delta \theta$ in Equation \eqref{eq:deltatheta} where $N$ represents the number of revolutions.
We graph the 0, 1, and 2 complete revolution solutions to Lambert's problem for this case with a time constraint of $59,952$s in Figure \ref{fig:MultipleRevolution}.
The specific energies for the 0, 1, and 2 complete revolution orbit transfers are -5.802 kJ/kg, -6.243 kJ/kg, and -9.958 kJ/kg respectively. Hence, the specific energy decreases as the number of revolutions increase of the transfer orbit for the specified transfer time.

\subsection{With $J_2$ Perturbations}
\begin{figure}
\begin{center}
\includegraphics[width=0.9\linewidth]{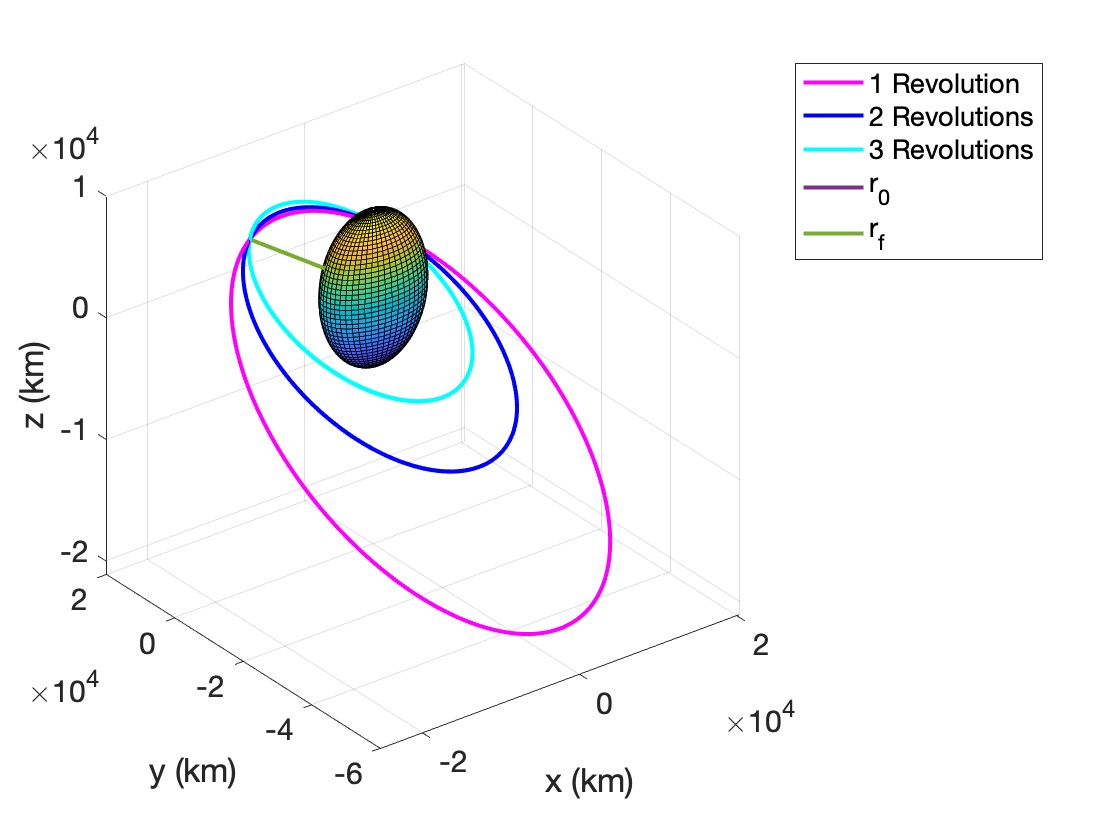}
\caption{Multiple Revolutions Orbit Trajectories with $J_2$ Effects for $t = 59,952$s}
\label{fig:MultipleRevolutionJ2}
\end{center}
\end{figure}

We take into account the oblateness of the Earth and find the multiple revolutions solutions within the time constraint of $t = 59,952$s. In this case, we plot the 1, 2, and 3 revolution solutions in Figure \ref{fig:MultipleRevolutionJ2}.
A benefit of using the KO methodology is that multiple revolutions solutions can easily be calculated by increasing the number $N$ without having to change the time of flight. 
This entails benefits because the zonal representation uses true anomaly as its independent variable.

\section{Conclusion}
This paper has demonstrated the application of Koopman operator theory to solve Lambert's problem for a satellite orbiting Earth both with and without taking into account $J_2$ perturbations. 
Approximate analytical solutions to both single and multiple revolutions Lambert's problem were found.
This paper builds upon previous results made by Arnas and Linares \cite{arnaslinares} to use operator theory to solve the zonal harmonics problem. A summary of the KO theory and appropriate variable transformations using state transition matrices was presented.

It was shown that solutions to Lambert's problem can be found for different initial conditions without the need to recalculate the KO matrix.
This provides a computational advantage of using the KO theory to solve Lambert's problem over traditional Lambert solvers.
It is important to note that the accuracy of the solution to Lambert's problem found using the KO decreases as the time increases and the number of revolutions increase upward of around $100$. This error can reach above 3 km in such situations. Further research should be made to compare the accuracy and computation speed of the KO method to traditional solvers of Lambert's problem. Also, different cost functions for solving Lambert's problem, such as minimizing Delta-V or minimizing transfer time, should be explored to find different optimal orbit trajectories.

\section*{Acknowledgments}
The authors wish to acknowledge the support of this work by the Defense Advanced Research Projects Agency (Grant N66001-20-1-4028). The content of the information does not necessarily reflect the position or the policy of the Government. No official endorsement should be inferred. Distribution statement A: Approved for public release; distribution is unlimited. The authors wish to acknowledge the support of this work by the MIT Portugal Program (MPP), Agmt Effective 6/1/18, and by the Air Force Office of Scientific Research (Award number R1302889).





\bibliographystyle{IEEEtran}
\bibliography{references.bib}
\end{document}